\documentclass[a4paper,oneside]{amsart}
\usepackage{amssymb}
\usepackage{amsthm}
\usepackage{amsmath,amscd}
\usepackage[mathscr]{euscript}
\usepackage[all]{xy}
\usepackage{float}
\usepackage[utf8]{inputenc}
\usepackage{lmodern}
\usepackage[T1]{fontenc}
\usepackage[textwidth=14cm,hcentering]{geometry}
\usepackage[colorlinks=true,linkcolor=red,citecolor=blue]{hyperref}
\usepackage{mathtools}
\usepackage{tikz}
\usepackage{pgf}
\usepackage{todonotes}
\usetikzlibrary{cd}
\usetikzlibrary{arrows,positioning}

\usepackage{hyperref}
\usepackage{cleveref}
\newtheorem{thm}{Theorem}[section]
\newtheorem*{thm*}{Theorem}

\newtheorem{rmq}[thm]{Remark}
\newtheorem{prop}[thm]{Proposition}
\newtheorem*{prop*}{Proposition}
\newtheorem{defn}[thm]{Definition}
\newtheorem{lemme}[thm]{Lemma}
\newtheorem{coroll}[thm]{Corollary}

\newtheorem{thmalpha}{Theorem} 
\Crefname{defn}{Definition}{Definitions}

\usepackage{mathtools}
\DeclarePairedDelimiter\ceil{\lceil}{\rceil}
\DeclarePairedDelimiter\floor{\lfloor}{\rfloor}

\newcommand{\gui}[1]{``#1''}
\newcommand{\xto}[1]{\xrightarrow{#1}}

\newcommand{\xinj}{\xhookrightarrow}

\newcommand{\inj}{\hookrightarrow}

\newcommand{\bb}[1]{\mathbb{#1}}
\newcommand{\mc}[1]{\mathcal{#1}}
\newcommand{\msf}[1]{\mathsf{#1}}
\newcommand{\mbf}[1]{\mathbf{#1}}

\newcommand{\bbZ}{\mathbb{Z}}

\newcommand{\vide}{\varnothing}

\newcommand{\ul}{\underline}

\newcommand{\bs}{\backslash}
\usepackage{mathtools}
\DeclareUnicodeCharacter{202F}{ }
\DeclareMathOperator{\Hom}{Hom}
\DeclareMathOperator{\Ext}{Ext}

\DeclareMathOperator{\colim}{colim}
\DeclareMathOperator{\holim}{holim}
\DeclareMathOperator{\hocolim}{hocolim}

\DeclareMathOperator{\conf}{Conf}
\DeclareMathOperator{\fiber}{fiber}
\DeclareMathOperator{\tfib}{tfib}
\DeclareMathOperator{\tcof}{tcof}
\DeclareMathOperator{\cofiber}{cofiber}
\DeclareMathOperator{\coker}{coker}

\usetikzlibrary{tikzmark}
\newcommand{\Map}{\text{Map}}
\newcommand{\FI}{\mathsf{FI}}
\newcommand{\FB}{\mathsf{FB}}
\usepackage{diagbox}
\begin{document}
\pagecolor{white}
\author[]{Nicolas Guès}
\address{Université Sorbonne Paris Nord, Laboratoire Analyse, Géométrie et Applications, CNRS (UMR 7539), 93430, Villetaneuse, France.}
\email{gues@math.univ-paris13.fr}
\thanks{This project is part of the PHC Pessoa 50834WE program. The author was supported by funding from Ecole Normale Supérieure and Ecole Doctorale Galilée (Université Sorbonne Paris Nord) during the writing of this project. }

\color{black}
	\title[linear ranges for homotopy groups of configuration spaces]{Linear stable ranges for integral homotopy groups of configuration spaces}
	\maketitle
  \begin{abstract}
  	We prove explicit linear stable ranges for the $\FI$-modules $\Hom(\pi_p \conf M, \bb Z)$ and $\Ext(\pi_p \conf M, \bb Z)$ with $\conf M$ being the configuration co$\FI$-space of a $d$-dimensional manifold with $d \geq 3$. The proof of this result uses a homotopy-theoretic approach to representation stability for $\FI$-modules. This allows us to derive representation stability results from homotopy-theoretical statements, in particular the generalized Blakers-Massey theorem. We also generalize to $\FI_G$-modules and orbit configuration spaces.
  	\end{abstract}
 \setcounter{tocdepth}{1}
 \tableofcontents

 \section*{Introduction}
 \subsection*{Representation stability and $\FI$-homology}
 The domain of representation stability deals with understanding the behaviour of coherent families of representations of the symmetric groups. Representation stability was introduced by Church and Farb in the foundational article \cite{church_representation_2013}, as a replacement of the notion of homological stability in this context. Many situations of this type arise from objects called $\FI$-modules, introduced in \cite{church_fimodules_2012}: these are functors from $\FI$ (the category of finite sets and injections) to an abelian category $\mc A$. Two important examples are:
 \begin{enumerate}
     \item{The homology of congruence subgroups: these are the kernels of $\mathrm{SL}_n(\bb Z) \to \mathrm{SL}_n(\bb F_p)$ for $p$ a prime, see \cite{church_fimodules_2014, putman_congruence_2015}.}
     \item{The cohomology of ordered configuration spaces of manifolds: to a finite set $S$, we associate the space $\mathrm{Emb}(S, M)$ where $M$ is a fixed manifold. If $S \inj T$ is an injection, precomposition defines a map $\mathrm{Emb}(T,M) \to \mathrm{Emb}(S,M)$. We get a functor $\FI^{op} \to \mathrm{Spaces}$ and thus the cohomology of configuration spaces naturally forms an $\FI$-module for each homological degree.}

 \end{enumerate}
We will write $V_S := V(S)$ for the evaluation of an $\FI$-module $V$ at the set $S$. We say that an $\FI$-module $V$ is finitely generated if there is a finite set $S$ of elements of $V$ such that no proper sub-$\FI$-module of $V$ contains $S$. An important feature of finite-generation is that an $\FI$-module over $\bb Q$ is finitely generated if and only if the underlying sequence of $S_n$-representations is representation stable in the sense of Church and Farb \cite{church_fimodules_2012}. Therefore, it is interesting to develop techniques to prove that a given $\FI$-module is finitely generated and control its generation degree. Towards this goal, a homology theory for $\FI$-modules has been developed in \cite{church_homology_2017}, which works over any base ring. For an $\FI$-module $V$, we define $\bb H_0 V$ as an $\FB$-module given on a finite set $T$ by

    $$( \bb H_0 V)_T  = \coker(\bigoplus_{S \subsetneq T} V_S \to V_T).$$
    
This gives a functor $\bb H_0 : \mbf{Fun}(\FI, \mc A) \to  \mbf{Fun}(\FB, \mc A)$. The $p$-th $\FI$-homology functor $\bb H_p$ is defined to be the $p$-th left derived functor of $ \bb H_0$.
For $V$ an $\FI$-module, call $$t_k V = \max  \{n \text{ such that } \bb H_k(V)_{\ul n} \neq 0\}.$$ The $\FI$-module $V$ is generated in finite degree (resp. presented in finite degree) if and only if $t_0 V < \infty$ (resp. $t_1 V < \infty$). An interesting consequence of finite generation/presentation is the following fact \cite[Theorem C]{church_homology_2017}. Take $N = \max{(t_0 V, t_1 V)}$. For any set $T$ of cardinality $n \geq N$, we have
$$ V_{T} \simeq \underset{S \subseteq T, |S| \leq N}{\colim} V_S.$$ 
In particular the functor is determined by its restriction to finite sets of cardinality $\leq N$.  \\
The main goal of this paper is to connect the theory of $\FI$-modules with the calculus of homotopy (co)cartesian cubes, which is widespread in functor calculus. We will show that many theorems in $\FI$-module theory can in fact be interpreted as connectivity properties of some cubical diagrams arising from the $\FI$-module structure. The connection comes from the following formula. For $V$ an $\FI$-module, the total $\FI$-hyperhomology $\bb H$ can be expressed as a total homotopy cofiber of a cube:
    $$\bb H(V)_T \simeq \cofiber ( \underset{ S \subsetneq T} {\hocolim} V_S \to V_T) =: \underset{S \subseteq T} \tcof V_S $$
One can extend the definition of $\bb H$ to $\FI$-objects in any pointed infinity-category with finite colimits and the formula above remains true.
One application of this homotopical approach is the following theorem.

\begin{thmalpha}\label{thmA}
	Let $M$ be a connected topological manifold of dimension $d \geq 3$. Then for $p\geq 2$, the homotopy groups $\pi_p(\conf M)$ can be given a co$\FI$-module structure, and the dual homotopy groups: $\Hom(\pi_p(\conf M), \bbZ) $ (non-torsion part) and $\Ext(\pi_p(\conf M), \bbZ) $ (torsion part) are $\FI$-modules presented in finite degree with the following bounds:
	$$ t_0 \Hom(\pi_p(\conf M), \bbZ) \leq 2\floor*{\frac{p-1}{d-2}}+ 2$$
	$$ t_1 \Hom(\pi_p(\conf M), \bbZ) \leq 2\floor*{\frac{p-1}{d-2}} + 3$$
		$$ t_0 \Ext(\pi_p(\conf M), \bbZ) \leq 2\floor*{\frac{p-1}{d-2}}+ 2$$
	$$ t_1 \Ext(\pi_p(\conf M), \bbZ) \leq 2\floor*{\frac{p-1}{d-2}} + 3$$
\end{thmalpha}
\begin{rmq}
	Note that the inclusion of co$\FI$-spaces $\conf_\bullet M \subset M^\bullet$ is $(d-1)$-connected by \cite[Theorem 2.2]{GGLG_2017}. Therefore, if $p < d-1$, then $\pi_p \conf_\bullet M \simeq (\pi_p M)^\bullet$. Hence the previous theorem gets improved :
\[
t_0 \Hom(\pi_p(\conf M),\bbZ) \le 1 \qquad
t_1 \Hom(\pi_p(\conf M),\bbZ) = -1
\]
\[
t_0 \Ext(\pi_p(\conf M),\bbZ) \le 1 \qquad
t_1 \Ext(\pi_p(\conf M),\bbZ) = -1
\]
	
\end{rmq}
Stability for (linear duals of) homotopy groups of configuration spaces over $\bb Z$ had been obtained in \cite[Theorem 1.3]{kupers_representation_2018} applying the Postnikov tower to the co$\FI$-space of configurations on $M$. However, explicit bounds were given only over $\bb Q$ using Quillen's approach to rational homotopy theory and the Eilenberg-Moore spectral sequence.
\\
Moreover, this derived approach gives improved stability bounds for the integral cohomology of configuration spaces, which are very similar to those recently obtained by Bahran \cite[Theorem G]{bahran_regularity_2024}.
\begin{thmalpha}\label{thmB}
	Let $M$ be a $u$-connected ($u\geq 0$) topological manifold of dimension $d \geq 3$. Then its integral cohomology groups are presented in finite degree with the following bounds.
	\begin{enumerate}
		\item{$t_0 H^p \conf M \leq \begin{cases} 
				2\floor*{\frac{2p+1}{d-1}} , & \text{if } u+1 \geq \frac{d-1}{2}, \\[6pt]
				2\floor*{\frac{p}{u+1}} & \text{otherwise.} \\[6pt]
			\end{cases}$}
		\item{$t_1 H^p \conf M \leq \begin{cases} 
				2\floor*{\frac{2p+1}{d-1}}+1, & \text{if } u+1 \geq \frac{d-1}{2}, \\[6pt]
				2\floor*{\frac{p}{u+1}}+1& \text{otherwise.} \\
			\end{cases}$}
	\end{enumerate}
\end{thmalpha}

The main motivation for adopting a derived approach to representation stability is that it is sometimes easier to get stability results for an $\FI$-complex than for its homology groups. If we know bounds for the values of $\mathbf t_k$ (the degrees of the $\FI$-hyperhomology, see \Cref{def:tk}) for a non-negatively graded $\FI$-complex $C$, the following theorem of Gan and Li allows to descend to the level of homology groups:
\begin{thm}\cite[Theorem 5]{gan_linear_2019}\label[thm]{thm:ganli}
	Let $W$ be a non-negatively graded  $\FI$-chain complex. Then
	$$ t_0(H_k(W)) \leq 2 \mbf t_k(W) + 1$$
	$$t_1(H_k(W)) \leq 2 \max\{\mbf t_k(W), \mbf t_{k+1}(W)\} + 2$$\end{thm}
Using this point of view, we will focus on an approach to improve representation stability ranges for configuration spaces developed by Miller and Wilson in \cite{miller_fihyperhomology_2019}. 
Their approach is the following: Using the \emph{puncture resolution} of configuration spaces of manifolds introduced by Randal-Williams in \cite{randal-williams_homological_2013}, they provide a bound for the $\FI$-hyperhomology of the cochain complex $C^*(\conf M, \bb Z)$ and then apply a variant of Gan and Li's theorem to get linear bounds for the $\FI$-homology of $H^*(\conf M)$. \\

 In particular we will extend this derived approach to the world of (co)$\FI$-spectra or (co)$\FI$-spaces. We study $\FI$-objects in a general $\infty$-category and show a \gui{Gan and Li-type} theorem (which bounds $t_\bullet \pi_n X$ in terms of bounds on $t_\bullet X$) when the category is stable. The proof adopts a different approach compared to previous related theorems, making use of the \gui{cardinality filtration} introduced in \Cref{FIgeneralities}, and gives linear bounds. A useful feature of this approach is that it also generalizes to co$\FI$-spectra, co$\FI$-spaces, and positively graded \emph{cochain} complexes. This level of generality will be necessary to deal with configuration co$\FI$-spaces. The previous approaches did not allow to deal with co$\FI$-objects nor nonconnective chain complexes, as we will explain later. \\

\subsection*{Outline}In \Cref{FIgeneralities}, we develop a general homotopy-theoretic approach to $\FI$-modules and along the way recover some known results in $\FI$-module theory with compact proofs. We introduce the \emph{cardinality filtration} of an $\FI$-object which we then apply in \Cref{GoingDownSection} to prove the  ``going-down`` theorem, relating the stability of an $\FI$-spectrum and the stability of its homotopy groups. We also prove a dual version for co$\FI$-spaces and co$\FI$-spectra. Finally we apply these general results in \Cref{Applications} to obtain a proof of \Cref{thmA} in the simply connected case and \Cref{thmB}. The core of the argument relies on the generalized Blakers-Massey theorem \cite[Theorem 2.5]{goodwillie_calculus_2003}  demonstrating its utility as a tool for deriving representation stability ranges for $\FI$-modules arising from a geometric context. Lastly, we show how this approach extends to the category $\FI_G$ and in particular to homotopy groups and cohomology of orbit configuration spaces. This allows to deal with the non-simply connected case in \Cref{thmA}.

\subsection*{Conventions}
For clarity, we make explicit our conventions of what are $n$-connected spaces and $n$-connected maps, following \cite[Definition 5.1.1]{munson_cubical_2015}. We will say that a space $X$ is $n$-connected if $\pi_i(X,x) = 0$ for every $0 \leq i \leq n$ and every basepoint $x$. A map $X \to Y$ is said to be $n$-connected if $\mathrm{hofiber}_y(f)$ is $(n-1)$-connected for every $y \in Y$. 
In particular the terminal map $X \to *$ is $(n+1)$-connected if and only if $X$ is $n$-connected.\\
In the rest of the paper, we will write $\ul n$ for the set $\{0,...,n-1\}$. Finally let $Q_n$ be the poset category of $\mc P(\ul n)$. A functor $Q_n \to \mc C$ is called an $n$-cube in $\mc C$.
For $X_\bullet$ an $n$-cube in a category $\mc C$ with a suitable notion of connectedness (such as spaces, or a stable category equipped with a $t$-structure). We will say that 
	\begin{enumerate}
		\item{$X_\bullet$ is $k$-cartesian if the natural map $X_\vide \to \holim_{\vide \neq S \subseteq \ul n} X_S$ is $k$-connected.}
		\item{$X_\bullet$ is $k$-cocartesian if the natural map $\hocolim_{S\subsetneq \ul n} X_S \to X_n$ is $k$-connected.\\}
	\end{enumerate}

\textbf{Acknowledgements.} I sincerely thank my supervisor, Geoffroy Horel, who suggested the idea of this project, for his invaluable feedback and insightful discussions.  I also thank Jeremy Miller, Aurélien Djament and Cihan Bahran for useful conversations and comments. I thank the anonymous referee for their careful proofreading and suggestions that improved some of the bounds and simplified parts of the arguments.
 \section{$\FI$-objects and $\FI$-homology}\label{FIgeneralities}
 \subsection{$\FI$-objects in an $\infty$-category}
 Let $\mc C$ be a cocomplete $\infty$-category.  An $\FI$-object in $\mc C$ is a functor $\FI \to \mc C$. These $\FI$-objects form an $\infty$-category $\FI \mc C$. In the same way we define $\FB \mc C$ the category of $\FB$-objets in $\mc C$. The inclusion functor $\FB \xto{i} \FI$ induces a pullback functor $i^*: \FI \mc C \to \FB \mc C$ which admits a left adjoint $M$ that will be described below. For $X$ an $\FB$-object, we will call $M(X)$ the \emph{free} $\FI$-object on $X$. In the $\FI$-module literature, free $\FI$-modules are sometimes called \emph{induced}.  
  \begin{prop}\label{freedescription}
     If $X$ is an $\FB$-object, let $M(X)$ be the $\FI$-module defined by:
     $$ M(X)_T = \bigoplus_{S \subseteq T} X_S$$
     if $f: T \to T'$ is an injection, we define $M(f)$ to send the summand $X_S$ to $X_{f(S)}$ in $M(X)_{T'}$. One can check easily that $M$ is functorial in $X$.  
     Then $M$ is a left adjoint to the forgetful functor $i^*$.
 \end{prop}
 \begin{proof}
     The proof in the abelian case found in \cite[Definition 2.2.2]{church_fimodules_2012} adapts verbatim.
 \end{proof}
 Note that $M$ really is a \emph{derived} left adjoint: in the classical abelian case, the (non-derived) left adjoint to $i^*$ is called $M$ too in the literature and happens to be exact. Therefore its total left derived functor is equal to itself and we can write $M$ in both cases without worry. We will write $M(m)$ for the representable $\FI$-space $\FI(\ul m, -)$. It is easy to see that $M(m)$ is equivalent to $M(\Sigma_m)$ where $\Sigma_m$ is seen as a discrete $\FB$-space concentrated in cardinality $m$. \\
 
If $\mc C$ is pointed, and $X$ is an $\FB$-object in $\mc C$, we can also construct an $\FI$-module $\alpha^* X$ which is the same as $X$ on objects and bijections, and maps non-bijective maps to the zero map. We will call $\alpha^*$ the \emph{extension-by-zero} functor. Note that $\alpha^*$ is not stricto sensu the pullback along a functor $\alpha : \FB \to \FI$. We will later construct a functor $\alpha : \FB_+ \to \FI_+$between the canonical enrichment over pointed sets of $\FB$ and $\FI$. This explains the terminology chosen here.

\subsection{The $\FI$-homology functor $\bb H$}
The extension-by-zero functor $\alpha^*: \FB \mc C \to \FI \mc C$ commutes with limits, and when $\mc C$ is presentable, it admits a left adjoint $\bb H:= \alpha_!: \FI \mc C \to \FB \mc C$ that we will call the $\FI$-homology functor. When $\mc C = \msf{S}_*$ (the $\infty$-category of spaces) or $\mc C= \msf{Sp}$ (the $\infty$-category of spectra) we will write $\bb H_p X$ for the $\FI$-module $\pi_p \bb H X$.
\subsection{Recovering usual $\FI$-homology}
Let $R$ be a ring. It is easy to check that the extension-by-zero functor $\alpha^*: \mathrm{Fun}(\FB ,R\mathrm{Mod}) \to \mathrm{Fun}(\FI ,R\mathrm{Mod})$ has a ($1$-categorical) left adjoint given by $$\bb (H_0V)_T = \coker(\bigoplus_{S \subsetneq T} V_S \to V_T).$$ Its total left derived functor is the $\FI$-hyperhomology $\mathbb H$ used in \cite{gan_linear_2019, miller_fihyperhomology_2019}. Therefore in the derived context, the derived left adjoint $ \bb H: \FI \mc D(R) \to \FB \mc D(R)$ to $\alpha^*$ recovers usual $\FI$-homology in the sense that $H_k \bb H V = \bb H_k V $ for an $\FI$-module $V$. 
\begin{defn}
	
	\label[defn]{def:tk}
	The \emph{degree} of an $\FI$-object $V$ (or $\FB$-object) is the maximum cardinality (potentially infinite) for which $V$ is non-zero.
	$$ \deg V = \max \{|S| \text{ such that } V_S \neq 0\} \in \bb N \cup \infty$$
	Define then 
	$$ \mathbf{t}_k V = \deg \bb H_k V.$$
	By convention $\deg V = -1$ for $V = 0$. 
\end{defn}
\subsection{Homotopical description of $\bb H$}
In this section we provide a homotopical description of $\bb H$ as the total cofiber of some cubical diagram arising from the $\FI$-module structure. \\
To begin, notice that the extension-by-zero functor $\alpha^*$ is not the pullback along a functor $\FI \to \FB$. If there was such $\alpha$, we could describe $\bb H$ as a left Kan extension of $\alpha$ and use the coend formula for Kan extensions. This problem can be solved by enriching $\FB$ and $\FI$ in pointed sets:
\begin{rmq}
    For $I$ a small category, write $I_+$ for its canonical enrichment over pointed sets: $\Hom_{I_+}(S, T)  = I(S,T)_+$. Then, for $\mc C$ a pointed $\infty$-category, there is an equivalence of $\infty$-categories
    $$\mbf{Fun}(I, \mc C) \simeq \mbf{Fun}_{*}(I_+, \mc C)$$
    where $\mc C$ on the right is seen as enriched over $\msf {S}_*$ and $\mbf{Fun}_*(-,-)$ is the space of pointed functors.  
\end{rmq}
Therefore, any $\FI$-object in a pointed $\infty$-category $\mc C$ can equivalently be seen as a pointed functor $\FI_+ \to \mc C$. Now the extension-by-zero functor $\alpha^*: \FB \mc C \to \FI C$ is simply the pullback along $\alpha: \FI_+ \to \FB_+$ which is the identity on objects, sends every bijection on itself and every non-bijective map to zero.\\

The next theorem is the main point of this section: it is a description of $\FI$-homology as a total cofiber. This allows to recover and understand in a more homotopical sense certain properties of $\FI$-homology, such as the long exact sequence of Gan \cite[Theorem 1]{gan_long_2016} or the chain complex computing $\FI$-homology \cite[Proposition 5.10]{church_homology_2017}. 
\begin{thm}
    \label{HomologyFormula}
    Let $V$ be an $\FI$-object in a pointed $\infty$-category $\mc C$. For any finite set $T$,
    $$(\bb HV)_T \simeq \cofiber ( \underset{ S \subsetneq T} {\hocolim} V_S \to V_T) =: \underset{S \subseteq T} \tcof V_S $$
\end{thm}
\begin{proof}

Using the coend formula for Kan extensions \cite{loregian_coend_2021}, we can describe the derived functor evaluated at $V$ as a coend (in the $\infty$-categorical sense):
$$(\bb HV)_T= \int^S \FB(S , T)_+ \overset{} \otimes V_S$$
which we will simply denote as a tensor product
$$ (\bb HV)_T = \FB(-, T)_+ \otimes_{\FI} V.$$
Now, consider the sequence of $\FB$-spaces:
$$\underset{S \subsetneq T} {\hocolim} \FI(- , S)_+ \xto{\iota} \FI(-,T)_+ \xto{\alpha} \FB(- ,  T)_+$$
We claim that it is a cofiber sequence. Let us check it pointwise and fix $U$ a finite set. If $|U| > |T|$, then the first term $\underset{S \subsetneq T} {\hocolim} \FI(U , S)_+ \simeq |N_\bullet(\mc P(T) \backslash \{T\} ) | \simeq *$ is a point,  as are the second and third terms. If $|U| = |T|$, then the first term is still a point, and the second and third terms are isomorphic to $S_n$ and $\alpha$ is the identity. Finally, if $|U| < |T|$, note that the functor
$$ \mc P(T) \backslash \{T\} \to \text{Set}$$
$$ S \mapsto FI(U,S)_+$$
can be written as a disjoint union of representables: 
$$S \mapsto  \{*\} \sqcup \bigsqcup_{i: U \inj T} \mathbf 1_{i(U) \subseteq S} $$ where $\mathbf 1_{i(U) \subseteq S}$ is either empty or a singleton whether $i(U) \subseteq S$. In other words $\mathbf 1_{i(U) \subseteq S}$ is representable by $i(U)$  and it is easy to see that its homotopy colimit is a point (equal to its strict colimit as a set). Therefore

$$\underset{S \subsetneq T} {\hocolim} \FI(U , S)_+ \simeq \underset{S \subsetneq T}{\colim} \FI(U,S)_+ \simeq \FI(U,T)_+.$$
And the last term is a point, so once again we get a cofiber sequence.
As derived coends commute with homotopy colimits in both variables we get a cofiber sequence of $\FB$-objects
$$\hocolim_{S\subsetneq T} (\FI(-,S)_+ \otimes V) \to \FI(-,T)_+\otimes V \to \FB(-,T)_+ \otimes V$$
which by the Yoneda lemma can be written
$$ \hocolim_{S \subsetneq T} V_S \to V_T \to (\bb H V)_T$$ and finishes the proof.
\end{proof}
We can then recover easily the following known result (see \cite[Theorem C]{church_homology_2017}) with a streamlined proof:
\begin{coroll}
	\label{corollim}
Let $V$ be an $\FI$-module, and $N = \max(t_0 V, t_1 V)$ which we assume to be finite. Then for any set $T$ with cardinality $n \geq N$, we have an isomorphism
$$
 \underset{S \subsetneq T, |S| \leq N} \colim{V_S} \xrightarrow{\sim} V_T.  
$$

Moreover $N$ is the smallest integer such that \ref{corollim} holds for any finite set $T$.
\end{coroll}
\begin{proof}
By \Cref{HomologyFormula}, we have the following fiber sequence of $\FI$-complexes
$$ \underset{ S \subsetneq T} {\hocolim} V_S \to V_T \to \bb H V (T)$$
We thus obtain an exact sequence in homology
$$0 \to H_1( \bb H V_T) \to H_0(\underset{ S \subsetneq T} {\hocolim} V_S) \to H_0(V_T) \to H_0(\bb HV_T)\to 0$$
which can be rewritten simply as
$$ 0 \to \bb H_1(V)(T) \to \underset{ S \subsetneq T} {\colim} V_S \to V_T \to \bb H_0(V)(T) \to 0$$
If $|T| \geq N$, then the extreme terms are zero if and only if the middle map is an isomorphism, and we obtain in this case
$$ \underset{ S \subsetneq T} {\colim} V_S \xto{\simeq} V_T.$$
We then obtain the expected result by induction. Indeed,
\begin{align*}
\underset{ S \subsetneq T} {\colim} V_S &= \underset{ U \subsetneq T} {\colim}\underset{ S \subsetneq U} {\colim} V_S\\
&= \underset{ U \subsetneq T} {\colim}\underset{ S \subsetneq U, |S| \leq N} {\colim} V_S \text{     by induction hypothesis} \\
&= \underset{ S \subsetneq T, |S|  \leq N} {\colim} V_S
\end{align*}
and the case $|T| = N$ is obvious.
\end{proof}
\begin{rmq}
    It is actually an unexpected feature that the value of $\bb H X$ for $X$ an $\FI$-object does not depend of the action of the symmetric groups on $X$ but only on the non-bijective maps, and that moreover one can compute $\bb H X$ at $T$ only from the maps arising from the inclusions $S \subseteq T$. This is a key property that we will use many times in this paper. This is closely linked to the fact that $\FI(n,m)$ for $n <m$ is the free $S_n$-set on the strictly increasing maps $\Delta_+(n,m)$ between $n$ and $m$.
\end{rmq}
Now we study how the free functor $M: \FB \mc C \to \FI \mc C$ interacts with $\FI$-homology.
\begin{prop}
    If $X$ is an $\FB$-object, then $\bb H M X \simeq X$.
\end{prop}
\begin{proof}
The composite of left adjoints $\bb H M: \FB \mc C \to \FB \mc C$ has a right adjoint $i^* \alpha^*  = (\alpha i)^*$ which is the pullback along $\alpha i: \FB_+ \to \FB_+$. This map is simply equal to the identity. Therefore $\bb H M \simeq id_{\FB \mc C}$.
\end{proof}
\subsection{co$\FI$ and co$\FB$-objects}
The configurations on a manifold naturally form a co$\FI$-space. If $\dim M \geq 3$ and $M$ is simply connected, its homotopy groups can be endowed with the structure of a co$\FI$-module \cite{kupers_representation_2018}. To study their stability, it will be necessary to dualize some of the previous definitions to co$\FI$-objects, which we do here. The situation is formally dual so the proofs are exactly the same.
\begin{defn}
    Let $i^*: \FI^{op} \mc C \to \FB^{op} \mc C$ be the pullback of the inclusion $\FB^{op} \inj \FI^{op}$. We define the cofree functor $M: \FB^{op} \mc C$ as the \emph{right} adjoint of $i^*$ when it exists. It is the case when $\mc C$ is the $\infty$-category of spaces or spectra. Equivalently this could be defined as the left adjoint of $i^*$ seen as a functor on $\FI$-objects in $\mc C^{op}$. 
\end{defn}
From now on an \gui{object} will either mean pointed space or spectrum.
\begin{lemme}
Let $X$ be a co$\FB$-object. Consider the cofree $\FI$-object $M(X)$. Then, for $T$ a finite set,
$$ M(X)_T \simeq \prod_{S \subsetneq T} X_S$$ with the obvious maps induced by injections of finite sets. 
\end{lemme}
\begin{proof}
    The proof is strictly dual to the $\FI$ case. 
\end{proof}
Dually there is an $\FI$-cohomology $\bb H$ defined for co$\FI$-objects in a pointed category:
\begin{defn}
    Let $\alpha^*$ be the extension-by-zero functor $\FB^{op} \mc C \to \FI^{op} \mc C$. We write $\bb H$ for its right adjoint (which exists if $\mc C$ is the $\infty$-category of pointed spaces or spectra). We will call $\bb H$ the $\FI$-cohomology functor.
\end{defn}
 Note that the $\FI$-cohomology of a co$\FI$-object $X$ is simply the $\FI$-homology of $X$ seen in the opposite category. We write the same way $\bb H$ (resp. $M$) for $\FI$-homology and $\FI$-cohomology (resp. free and cofree object). This will not be ambiguous in most cases because they apply to different kinds of objects (even if $\FB \simeq \FB^{op}$ it will be clear from context whether we are working with $\FI$ or co$\FI$-modules).
\begin{lemme}
Dually to the $\FI$ case, we have:
\begin{enumerate}
    \item{for a co$\FI$-object $X$, $$ \bb H X (T) = \tfib_{S \subseteq T} X_S$$}
    \item{for a co$\FB$-object $X$, $\bb H M X = X$}
\end{enumerate}
\end{lemme}
\begin{rmq}
	Note that all of this makes sense only for pointed objects.  For co$\FI$-unbased spaces, the $\FI$-cohomology $\bb HX$ does not exist as a co$\FB$-object, because without a coherent choice of basepoint one cannot do the previous construction. However, for a set $S$, one could choose a basepoint in $X_S$ and propagate it to all the subsets $T \subset S$. Then the total fiber
	$$ \tfib_{T\subseteq S} X_T$$ makes sense. However, the collection of all the total fibers is not naturally a co$\FB$-module since the basepoints are not coherent with the action of the symmetric groups. In any case we can still define
	the $k$-th generation degree $\mbf{t}_k X$ by
	$$ \mbf t_k X = \max {\{n \in \bb N | \pi_k \tfib_{T\subseteq \ul n} X_T \neq 0 \text{ for some choice of basepoint } x_T \in X_T \}}$$
\end{rmq}
\subsection{Shift functor and $\FI$-derivative}
Taking disjoint union with a point defines an endofunctor of $\FI$. This gives rise by precomposition to the \emph{shift}  functor $S: \FI \mc C \to \FI \mc C $ where $\mc C$ is any category \cite{church_fimodules_2014}. There is a natural transformation $id \to S$. We can define the $\FI$-derivative functor $D$ to be the cofiber of $id \to S$. In \cite{gan_long_2016} a long exact sequence in $\FI$-homology relating the homology of an $\FI$-module and its shift is derived: 
\begin{thm*}\cite[Theorem 1]{gan_long_2016}
    Let $X$ be an $\FI$-module. There is a long exact sequence
    $$ ... \to \bb H_a X \to \bb H_a SX \to S \bb H_a X \to ...$$
\end{thm*}
We show below that it can be generalized to $\FI$-spectra or $\FI$-spaces: 
\begin{thm}\label{homologyles}
    Let $X$ be an $\FI$-spectrum, $SX$ the shift defined as above. Then there is an exact sequence of $\FB$-spectra (resp. a cofiber sequence in $\FB$-spaces):
    $$ \bb H X \to \bb H S X \to S \bb H X$$
    Which gives rise to a long exact sequence of homotopy groups
    $$ ... \to \bb H_a X \to \bb H_a SX \to S \bb H_a X \to ...$$
    Dually, if $X$ is a co$\FI$ spectrum or a \emph{pointed} co$\FI$-space, There is a fiber sequence 
    $$  S \bb H X \to \bb H S X \to \bb H X $$
    which gives rise to a long exact sequence of homotopy groups 
    $$ ... \to S \bb H_a X \to \bb H_a SX \to \bb H_a X \to ...$$

\end{thm}
In particular we recover Gan's result when applying it to the Eilenberg-Maclane $\FI$-spectrum $X = H V$ when $V$ is an $\FI$-module.
\begin{proof}
As $DX$ is the homotopy cofiber of $X \to SX$ and $\bb H$ commutes with cofiber sequences, it suffices to show that $\bb H DX = S \bb H X$. Fix an integer $n$. Notice that $\bb H X_{\ul {n+1}}$ is the total cofiber of an $(n+1)$-cube $\mc X$, that we will write as a map of $n$-cubes
$$ \mc Y \to \mc Z$$
where $\mc Y$ is simply the cube of $X_{I}$ for $I \subseteq\ul{n}$, whereas $\mc Z$ is the cube indexed by the sets $$ \{1\} \sqcup  S \subseteq\ul {n+1} $$ for $S \subseteq\{2,...,n\}$. Therefore, by proposition 5.9.3 in \cite{munson_cubical_2015} we get a cofiber sequence
$$ \tcof \mc Y \to \tcof \mc Z \to \tcof \mc X$$
that we can simply identify with:
$$ \bb H X_{\ul n} \to \bb H S X_{\ul n} \to \bb H X_{\ul {n+1}} = S \bb H X_{\ul n}$$
To conclude note that all these identifications and cofiber sequences used are natural with respect to bijections of finite sets and therefore lift to cofiber sequences of $\FB$-objects.
The dual proof works the same way.
\end{proof}
\subsection{A complex computing $\FI$-homology}
In \cite[Definition 2.19]{church_fimodules_2014} and \cite[Definition 5.9]{church_homology_2017}, two complexes were introduced to compute the $\FI$-homology of an $\FI$-module $V$, respectively $\overset{\sim}{S}_{-\bullet} V$ and $V \otimes_{FI} C_\bullet$ and are shown to be the same in \cite[section 5.1]{church_homology_2017}. In this section, we demonstrate how this complex can be recovered through a concise proof leveraging cubical calculus and a spectral sequence coming from \cite{munson_cubical_2015}.

Proposition 9.6.14 from \cite{munson_cubical_2015} gives a spectral sequence computing the homology of the total cofiber of a cube of spaces ; because the category of spectra is stable, the argument works for the the homotopy groups of a cube of spectra with the same proof and we get:
\begin{lemme}
    Let $X$ be a $n$-cube of spectra. There is a spectral sequence computing the homotopy groups of the total cofiber given by:
    $$E^1_{p,q}  = \bigoplus_{S \subseteq\ul n, |S| = n-p} \pi_q(X_S) \implies \pi_{p} \tcof{X} $$
with the first differential
\[
d^1:
\bigoplus_{S \subseteq P(n), \, |S| = n - p } \pi_q(X_S) \longrightarrow
\bigoplus_{T \subseteq P(n), \, |T| = n - p+1} \pi_q(X_T)
\]
whose restriction to the summand indexed by  $S$  is
\[
d^1|_{\pi_q(X_S)} =
\sum_{i \notin S} (-1)^{|S| \cdot i} X(S \to S \cup \{i\})
\]

\end{lemme}

Now let $X$ be an $\FI$-spectrum. Seeing $\FI$-homology as the total cofiber of the cubes arising from the $\FI$-object structure, we get a spectral sequence 
$$ E^1_{p,q} = \bigoplus_{S \subseteq\ul n, |S| = n-p} \pi_q(X_S) \implies \bb H_p X(n). $$
In the particular case when $X$ is concentrated in degree $0$ (i.e. if it comes from an $\FI$-module $V$ seen as an $H\bbZ
$-module spectrum), $E^1_{p,q}$ is simply a chain complex, because the terms are concentrated on the line $q= 0$ and thus the spectral sequence collapses at $E^2$. We get that for an $\FI$-module $V$, the $\FI$-homology $\bb H_\bullet V$ can be computed as the homology of the chain complex $S_\bullet$ with $$S_p = E^1_{p,0} =  \bigoplus_{S \subseteq\ul n, |S| = n-p} V_S.$$

A careful examination of the differential shows it is exactly the one from \cite{church_fimodules_2014} and \cite{church_homology_2017}.
\subsection{The cardinality filtration}
In this section we introduce a tool for the proof of the \gui{going-down} theorem \Cref{goingdown}, the \emph{cardinality filtration}. It realizes a canonical dévissage of an $\FI$-object $X$ by free objects. This filtration also appears in Kaya Arro's recent paper \cite[Definition 18]{arro_ficalculus_2023} on $\FI$-calculus.
\begin{defn}
    For an integer $k$, call $\FI_{\leq k}$ (resp $\FB_{\leq k}$) the restriction of $\FI$ (resp $\FB$) to sets smaller than $k$. The tower of inclusions
    $$ \FI_{\leq 0} \inj \cdots \inj \FI_{\leq k-1} \xinj{} \FI_{\leq k} \inj \cdots \inj \FI$$ yields a tower of restriction functors $\FI \mc C \to \cdots \to \FI_{\leq k} \mc C \xto{} \FI_{\leq {k-1}} \mc C \to \cdots \FI_{\leq 0 } \mc C $. 
\end{defn}
\begin{defn}
    For an integer $k$, call $(i_k)_!$ the (derived) left adjoint of the restriction functor $\FI \mc C \xto{(i_k)^*} \FI_{\leq k} C$. We will call $F_k$ the composite $\FI \mc C \xto{(i_k)^*} \FI_{\leq k} \mc C \xto{(i_k)_!} \FI \mc C$. It is the (derived) counit of the adjunction and we have by standard categorical nonsense a sequence of natural transformations $F_0 \to \cdots \to F_k \to F_{k+1} \to \cdots \to id_{\FI \mc C}$.
\end{defn}
\begin{defn}\label[defn]{cardinalityfiltration}
For $X$ an $\FI$-object, We call the tower 
$$F_0 X \to \cdots \to F_k X\to \cdots \to X$$
the cardinality filtration of $X$.
\end{defn}
Note that the cardinality filtration is always exhaustive because it is pointwise stationary: $(F_kX)_T \simeq X_T$ as soon as $|T| \leq k$. 
This filtration can be thought as some flavour of a $0$-dimensional Goodwillie-Weiss tower for presheafs over the category $Disk^0$ of finite disjoint unions of $0$-dimensional discs and embeddings (that is, finite sets and injections). \\

Now let us describe the layers of this tower.

\begin{lemme}\label[lemme]{lemma2.6}
    Let $X$ be a $\FI$-object in a pointed $\infty$-category $\mc C$     and $k$ an integer. Then:
    \begin{enumerate}
        \item{For $j \leq k$, $(\bb H F_k X)_{\ul j} \simeq (\bb H X)_{\ul j}$}
        \item{For $j > k$, $(\bb H F_k X)_{\ul j} = 0$.}
    \end{enumerate}
\end{lemme}
\begin{proof}
    Call $j_k$ the inclusion $\FB_{\leq k} \inj \FB$. There is an obviously commutative diagram:
    \[\begin{tikzcd}
	{\FI_{\leq k,+}} & {\FB_{\leq k, +}} \\
	{\FI_{+}} & {\FB_{ +}}
	\arrow["{i_k}", from=1-2, to=2-2]
	\arrow["{j_k}", from=1-1, to=2-1]
	\arrow["\alpha", from=1-1, to=1-2]
	\arrow["\alpha"', from=2-1, to=2-2]
\end{tikzcd}\]
This implies that the corresponding diagram on functor categories commutes; By standard category theory facts, this implies that their left adjoints commute and we get $ \bb H \circ (i_k){}_! \simeq (j_k){}_! \circ \bb H$, with a slight abuse of notation for $\bb H$. But $(j_k){}_!: \FB_{\leq k,+}\mc C \to \FB_+ \mc C$ is just an extension by zero, because there are no morphisms between non-isomorphic sets in $\FB$ and $\FB_{\leq k} \inj \FB$ is a fully faithful inclusion. Finally we note that $\bb H$ commutes with $(i_k)^*$ in the sense that the $\FI_{\leq k}$-homology of the restriction is the restriction of the $\FI$-homology. Putting this together we get:
$$ \bb H \circ F_k (X) = \bb H \circ (i_k){}_! \circ (i_k)^* X = (j_k){}_! \circ \bb H \circ (i_k)^* X  = (j_k){}_{!} \circ (j_k)^* \circ \bb H X $$ which is exactly $(\bb H X)_S$ for $|S| \leq k$ and $0$ otherwise.
\end{proof}
We can finally describe the layers of the cardinality filtration: 
\begin{thm}\label{layersfiltration}
    Let $X$ be an $\FI$-object. We have an equivalence of $\FI$-objects:
    $$ \cofiber (F_{k-1} X \to F_{k} X) \simeq M(\bb H X(\ul {k}))$$
    where $\bb H(X)(\ul {k})$ is seen as an $\FB$-module concentrated in cardinality $k$ (we take the convention $F_{-1} X = 0$). 
\end{thm}
\begin{proof}
Write $L_k X$ for the cofiber of $F_{k-1}X \to F_k$ and $\ell_i$ for the inclusion $\FI_{\leq i-1} \subset \FI_{\leq i}$. We have obviously $i_{k} \ell_k  = i_{k-1}$ and therefore 
$$(i_k)_! (i_k)^* = (i_k)_! (\ell_k)^* (i_k)^*$$
One can check that the natural map
$$ F_{k-1} X  = (i_{k-1})_! (i_{k-1})^* X \to (i_k)_! (i_k)^* X = F_{k} X$$
is simply the image of the counit
$$ (\ell_k)_! (\ell_k)^* \to id_{\FI_{\leq k}}$$
by the functor $(i_k)_!$, evaluated on the restriction $(i_k)^* X$. Consequently, we have that 
$$ L_k X \simeq \cofiber{\bigg( (i_{k+1})_! \bigg( (\ell_k)_! (\ell_k)^* (i_k)^* X \to (i_k)^* X \bigg) \bigg)}$$
as $(i_{k+1})_!$ is a left adjoint it commutes with homotopy cofibers, so we have
$$\cofiber(F_{k-1} X \to F_k X)  \simeq (i_{k+1})_! \cofiber{\bigg( (\ell_k)_! (\ell_k)^* (i_k)^* X \to (i_k)^* X \bigg)}$$
The cofiber on the right-hand side is a $\FI_{\leq k}$-object concentrated in cardinality $k$. Therefore, this means it comes from a $\Sigma_k$-object $Y$ (and extended by zero in cardinality different than $k$). From the commutative diagram
\[\begin{tikzcd}
	& {\FI_{\leq k}} \\
	{\Sigma_k} && \FI \\
	& \FB
	\arrow[from=1-2, to=2-3]
	\arrow[from=2-1, to=1-2]
	\arrow[from=2-1, to=3-2]
	\arrow[from=3-2, to=2-3]
\end{tikzcd}\]
we deduce a commutative diagram of left Kan extensions:
\[\begin{tikzcd}
	& \mbf{Fun}(\FI_{\leq k}, \mc C) \\
	{\mbf{Fun}(\Sigma_k, \mc C)} && {\mbf{Fun}(\FI, \mc C)} \\
	& {\mbf{Fun}(\FB, \mc C)}
	\arrow["{(i_k)_!}", from=1-2, to=2-3]
	\arrow["{\text{ext}}", from=2-1, to=1-2]
	\arrow["{\text{ext}}"', from=2-1, to=3-2]
	\arrow["M"', from=3-2, to=2-3]
\end{tikzcd}\]
Where the $\mathrm{ext}$s are the natural extension-by-zero functors. 
Consequently, we get that 
$$ L_k X \simeq M(Y)$$
So the layer $L_k X$ is a free $\FI$-object on a $\Sigma_k$-object. It remains to identify $Y$. This is easy since $\bb H L_k X = \bb H MY = Y$ and by \Cref{lemma2.6} $\bb H L_n X$ is equal to $(\bb H X )_{\ul k}$ in cardinal $k$ and zero elsewhere. 
\end{proof}
\section{Stability of a (co)$\FI$-object versus stability of its homotopy groups}\label{GoingDownSection}
Linking the stability bounds of an $\FI$-chain complex and the stability of its homology  proves useful: this is the approach used in \cite{miller_fihyperhomology_2019}, \cite{bahran_regularity_2024} to obtain improved stability bounds for ordered configuration spaces of manifolds, using Gan and Li's result. 
The goal of this section is to generalize this result to (co)$\FI$-spectra and co$\FI$-spaces. Our main result is that we obtain linear bounds in $k$ for the generation and presentation degrees of $\pi_k X$ if $X$ is a stable $\FI$-spectrum with generation degrees $\bf t_k$ increasing linearly, providing a generalization of Gan and Li's theorem. We will call that kind of result a \gui{going-down} theorem: from the higher notion of stability we recover stability at the linear level of the homotopy/homology groups. This kind of results has some history:
\begin{enumerate}
	\item{The first bounds obtained in \cite[Claim 5.14]{church_homology_2017} gave an  exponential bound (in $p$) on $t_0 H_p C$ and $t_1 H_p C$ when $C$ is a positively graded $\FI$-chain complex with $\mbf t_p C$ controlled by a affine function of $p$.}
	\item{This bound was widely improved in \cite[Theorem 5.1]{church_linear_2018}, yielding a quadratic bound, using the properties of the \emph{local degree} and the \emph{stable degree} for $\FI$-modules. }
	\item{In \cite[Theorem 5]{gan_linear_2019} a linear bound was obtained with different techniques, still for positively graded chain complexes and improved in \cite[Theorem 2.20]{bahran_regularity_2024}. This is the best ``going-down`` result known so far.}
\end{enumerate}
Our strategy of proof generalizes also to the dual setting of co$\FI$-spectra and co$\FI$-spaces. If $X$ is a co$\FI$-space/spectrum, recall that its $\FI$-cohomology $\bb H X$ is defined to be the collection of the total fibers of the cubes arising from the co$\FI$ structure. In this case, we will provide bounds for the stability of the dual homotopy groups of $X$ (which are $\FI$-modules) in term of the slope controlling the connectivity of $\bb H X$. It is necessary to deal with this dual setting, because the configuration spaces of a manifold form a co$\FI$-space. Note that a going-down theorem is this dual setting is not a formal consequence of the going-down for connective $\FI$-complexes: with an approach using the hyperhomology spectral sequence, there will be concrete convergence problems in trying to apply an argument like in \cite{church_homology_2017} or \cite{church_linear_2018}. 

Let us prove first an easy converse of the \gui{going-down} theorem. We will not use this result in the rest of the paper, but believe it is relevant to include it as a potiental tool within this homotopical $\FI$ machinery.
\newcommand{\hmax}{h^{\max}}
\begin{prop}[``Going-up``]
	Let $X$ be a globally bounded-below $\FI$-spectrum. Then $X$ has degrees
 $$ \mbf t_k X \leq  \underset{j \in \bb Z}\max \; t_{k-j} \pi_j X$$
\end{prop}
\begin{proof}
    We can assume that $X$ is a connective $\FI$-spectrum by suspending it enough times (suspension commutes with $\bb H$). Let $\tau_{\leq n}$ be the $n$-th truncation functor . There are natural transformations $\tau_{\leq n} \to \tau_{\leq {n-1}}$. For any $X \in \FI \mbf{Sp}$, we have 
    $$ X = \lim \tau_{\leq n} X$$
    and an exact sequence
    $$ \Sigma^n \pi_n X \to \tau_{\leq n} X \to \tau_{\leq {n-1}} X$$
    Since $\bb H$ is exact, we obtain a long exact sequence of homotopy groups and in particular an exact sequence 
    $$ \bb H_k \Sigma^n \pi_n X \to \bb H_k \tau_{\leq n} X \to \bb H_k \tau_{\leq {n-1}}$$
    It follows that 
    $$\mbf t_k (\tau_{\leq n} X) \leq \max (\mbf t_k(\tau_{\leq {n-1}} X), t_{k-n} \pi_n X)$$
    and thus by induction,
    $$ \mbf t_k(\tau_{\leq n} X) \leq \underset{j \leq n}\max \; t_{k-j} \pi_j X$$
    Finally, finite colimits commute with filtered limits in $\mbf{Sp}$ and thus $\bb H$ commutes with filtered limits. It follows that: 
    $$ \mbf t_k X \leq  \underset{j \in \bb N}\max \; t_{k-j} \pi_j X$$
\end{proof} 
\begin{thm}[Going-down theorem, $\FI$ version] \label{goingdown}
Let $X$ be an $\FI$-spectrum. Suppose that $\mathbf{t}_p:= \deg \bb H_p X < \infty $ for all $p$. Then for every $p \in \bb Z$, the $\FI$-module $\pi_p X$ is zero if $\mbf t_p =-1$ and is presented in finite degrees with the following bounds if $\mbf t_p \geq 0$ : 
\begin{enumerate}
    \item{    $$t_0 \pi_p X \leq  \max\{2 \mbf t_p, 2\mbf t_{p+1}-1, 2 \mbf t_{p+2} - 2\}$$
    	$$t_1 \pi_p X \leq  \max\{2 \mbf t_p+1, 2\mbf t_{p+1}, 2 \mbf t_{p+2} - 1 \}$$}
    \item{In the special case when $\mbf{t}_p \leq f(p)$ with $f$ weakly decreasing (the $\bb H X_S$ are highly \emph{coconnected}), this simplifies to : 
    $$ t_0 \pi_{p} X\leq 2f(p)$$
    $$ t_1 \pi_{p} X\leq  2f(p)+1$$  }
    \end{enumerate}
\end{thm}

\begin{thm}[Going-down, co$\FI$ versions\label{cogoingdown}]
Let $X$ be either: \begin{enumerate}
	 \item a connected pointed co$\FI$-space with $\pi_1(X_S)$ abelian for each $S$ or
	 \item a co$\FI$-spectrum.
	 \end{enumerate} 
	 Fix an injective abelian group $A$. Suppose $\mbf t_p = \deg \pi_p \bb HX < \infty$ for all k. Then for every $p$, the dual $\bb D \pi_p X:= \Hom(\pi_p X, A)$ is zero if $\mbf t_p = -1$, and is presented in finite degrees if $\mbf t_p \geq 0$ with 
\begin{enumerate}
    \item{$t_0 \bb D\pi_p X \leq  \max \{2\mbf t_p, 2\mbf t_{p-1} -1, 2\mbf t_{p-2}-1\}$}
      \item{$t_1 \bb D \pi_p X \leq  \max \{2\mbf t_p+1, 2\mbf t_{p-1}-1, 2\mbf t_{p-2}-2\}$}
\end{enumerate}
Additionally, if we have $\mbf t_p \leq f(p)$ with $f$ weakly increasing then we get 
        \begin{enumerate}
        \item{$t_0 \bb D \pi_{p} X \leq  2f(p)$}
        \item{$t_1 \bb D \pi_{p} X \leq  2f(p)+1$}
    \end{enumerate}
Moreover, when $X$ is a co$\FI$-space, the assumption of having a basepoint can be dropped if we assume that $X_S$ as well as $\bb HX_S$ are simply connected for all $S$.  
\end{thm}
To prove \Cref{goingdown} and \Cref{cogoingdown} we will use the spectral sequence arising from the cardinality filtration, so we need to control how the bounds $t_0, t_1$ for $\FI$-modules over $\bb Z$ behave through spectral sequences. This is the main point of \cite{church_linear_2018}. Their strategy is to compare $t_0$ and $t_1$ to other invariants $\delta$ (the stable degree) and $\hmax$ (the local degree).We will not recall their precise definition as the only properties needed here are the following facts. 
\begin{lemme}\cite[Proposition 3.1, Corollary 2.13]{church_linear_2018}\label[lemme]{localmaxbound}
We have the following relations between $\delta, h^{\max}, t_0, t_1$:
    \begin{enumerate}
	    \item{$\delta \leq t_0$}
	    \item{$h^{\max} \leq t_0+ \max(t_0,t_1) -1$.}
	    \item{If $V$ is a free $\FI$-module, then $\hmax V = -1$.}
	    \end{enumerate}
\end{lemme}
We will use the following recent improvement by Bahran of \cite[Proposition 4.1]{church_linear_2018}. For $V$ an $\FI$-module presented in finite degrees, define $\mathrm{reg} V = \sup_p t_p V - p$ which is finite by \cite{church_homology_2017}.  In particular $\mathrm{reg} V \geq t_1 V - 1$. ·
\begin{thm}[{\cite[Particular case of Corollary C]{Bah25}}]\label[theorem]{bahranspectral}
	Let $(U^r_k : r \ge r_0)$ be a homological single-graded spectral sequence of $\FI$-modules and $(\varepsilon_k : k \in \mathbb{Z})$ a $\mathbb{Z}$-indexed sequence of integers $\ge 1$ such that at the inital page $r_0$, for every $k \in \mathbb{Z}$ the $\FI$-module $U^{r_0}_k$ is presented in finite degrees with
	\begin{itemize}
		\item $\delta(U^{r_0}_k) \le \varepsilon_k$,
		\item $h^{\max}(U^{r_0}_k) = -1$.
	\end{itemize}
	Then for every $r \ge r_0$ and $k \in \mathbb{Z}$, the $\FI$-module $U^r_k$ at page $r$ satisfies the following:
	\begin{enumerate}
		\item $\delta(U^r_k) \le \varepsilon_k$.
		\item $h^{\max}(U^r_k) \le \max\{-1, 2\varepsilon_k - 2, 2\varepsilon_{k+1} - 2, 2\varepsilon_{k+2}-2 \}$.
		\item Both $t_0(U^r_k)$ and $\text{reg}(U^r_k)$ (and hence $t_1(U^r_k) - 1$) are bounded above by
		$$
		\begin{cases}
			\max\{-1, 2\varepsilon_{k+1} - 2, 2\varepsilon_{k+2} - 2\} & \text{if } \varepsilon_k = -1,\\
			\max\{2\varepsilon_k, 2\varepsilon_{k+1} - 1, 2\varepsilon_{k+2} - 2\} & \text{if } \varepsilon_k \ge 0.
		\end{cases}
		$$
	\end{enumerate}

\end{thm}

We will also need to control $t_0$ and $t_1$ in terms of $\hmax$ and $\delta$. The best bound known is the following. I thank Cihan Bahran for pointing out the existence of this bound, which leads to a quantitative enhancement of the theorem compared to the first version of this paper. 
\begin{lemme}\cite[Corollary 2.13, Theorem 2.6]{bahran_regularity_2024}\label[lemme]{bahranlemma}
			Let \( V \) be an $\FI$-module presented in finite degrees.  Then
			
			\[t_0(V) \leq 
			\begin{cases} 
				\delta & \text{if } \hmax = -1, \\ 
				\hmax & \text{if } \delta = -1 \text{ and } \hmax \geq 0, \\ 
				\hmax + 1 & \text{if } 0 \leq \delta \leq \lceil \hmax/2 \rceil \text{ and } \hmax \geq 0, \\ 
				\delta + \lfloor \hmax/2 \rfloor + 1 & \text{if } \delta > \lceil \hmax/2 \rceil \text{ and } \hmax \geq 0, 
			\end{cases}\]
			
			and
			
			\[t_1(V) \leq 
			\begin{cases} 
				-1 & \text{if } \hmax = -1, \\ 
				\hmax + 1 & \text{if } \delta = -1 \text{ and } \hmax \geq 0, \\ 
				\hmax + 2 & \text{if } 0 \leq \delta \leq \lceil \hmax/2 \rceil \text{ and } \hmax \geq 0, \\ 
				\delta + \lfloor \hmax/2 \rfloor + 2 & \text{if } \delta > \lceil \hmax/2 \rceil \text{ and } \hmax \geq 0. 
			\end{cases}\]

\end{lemme}

\begin{lemme}\cite[Proposition 3.3]{church_linear_2018}\label[lemme]{localmaxexact}
Let $A \to B$ be a map between $\FI$-modules presented in finite degrees. Then,
\begin{enumerate}
    \item{$\delta(\ker f) \leq \delta(A)$}
    \item{$\delta(\coker f) \leq \delta(B)$}
    \item{$h^{\max} (\ker f) \leq \max(2\delta A - 2, h^{\max} A, h^{\max} B)$}
    \item{$h ^{\max} (\coker f) \leq \max(2\delta A - 2, h^{\max} A, h^{\max} B)$}
\end{enumerate}    
\end{lemme}
\begin{proof}[Proof of \Cref{goingdown}]
    Let us study the spectral sequence associated to the cardinality filtration $(F_k X)$. We take as indexing convention that
    $$E^1_{p,q} = \pi_p(C_q X) = \pi_p(M \bb H (X)_ {\ul{q}}). $$
    The differential $d_r$ at the $r$-th page with this indexing convention is of bidegree $(-1, -r)$.  
    Note first that the $E^1_{p,q}$ are zero for $q < 0$. Moreover, as $C_q X = M \bb H (X)_{\ul{q}}
    $ is a free $\FI$-object, we have
    $$ \pi_p (C_q X)_S = \pi_p \left( \bigoplus_{T \subseteq S, |T| = q} \bb H X_T \right) =\bigoplus_{T \subseteq S, |T| = q} \bb \pi_p  \bb H X_T= M(\pi_p \bb H X_{\ul q}) $$
    Thus the $E^1$ page contains only free $\FI$-modules. For $q > \mbf t_p(X)$, $\bb \pi_p \bb H X _{\ul q}= 0$ by definition and thus $E^1_{p,q} = 0$ as soon as $q > \mbf t_p(X)$. Therefore each column is finitely supported. To sum up, the $E^1$ page looks like this (the red dots are nonzero terms).
\pgfmathsetseed{3345}
	
	\begin{figure}[H]
		\centering
		\begin{tikzpicture}[scale=0.8]
			\draw[->] (-5,0) -- (5,0) node[right] {$p$};
			\draw[->] (0,-0.5) -- (0,5.5) node[above] {$q$};
			\draw[->, blue] (2,3) -- (1,2) node[midway,above] {$d_1$};
			\draw[->, blue!50!cyan] (2,3) -- (1,1) node[midway,above] {$d_2$};
			\draw[->, blue!20!cyan] (2,3) -- (1,0) node[midway,above] {$d_3$};
			\foreach \x in {-5,...,5} { 
				\pgfmathtruncatemacro{\r}{random(1,4	)} 
				\foreach \y in {0,...,5} {
					\ifnum\y>\r
					\fill[gray] (\x,\y) circle (1pt); 
					\else
					\fill[red] (\x,\y) circle (2pt); 
					\fi
				}
			}
		\end{tikzpicture}
	\end{figure}
	    For all $r$, set $V_r^p = \bigoplus_{q=0}^{+\infty} E_r^{p,q}=\bigoplus_{q=0}^{\mbf{t}_p(X)} E_r^{p,q}$. By \Cref{localmaxbound}, the collection $(V^p_r)$ is a single-graded spectral sequence satisfying the hypothesis of \Cref{bahranspectral}. In particular we get $$t_0 \pi_p X \leq \max \{2 \mbf t_p, 2\mbf t_{p+1} -1, 2\mbf t_{p+2} - 2\}$$
	     $$t_1 \pi_p X \leq \max \{2 \mbf t_p+1, 2\mbf t_{p+1}, 2\mbf t_{p+2} - 1\}$$
Suppose now that $\mbf t_p \leq f(p)$ with $f$ weakly decreasing. The general case simplifies then to 
$$ t_0 \pi_p X \leq 2f(p)$$
$$ t_1 \pi_p X \leq 2f(p)+1$$
%
%

\end{proof}
\begin{rmq}\label[rmq]{rmq2.8}
	Notice that applying \Cref{bahranspectral} to the cardinality filtration spectral sequence also yields bounds for $\delta$ and $\hmax$ of $\pi_p X$ : 
	\begin{enumerate}
		\item $\delta \pi_p X \leq \mbf t_p,$
		\item $\hmax \pi_p X \leq \max\{-1, 2\mbf t_p-2, 2\mbf t_{p+1}-2, 2\mbf t_{p+2}-2\}.$
	\end{enumerate}
	Dually, if $X$ is a co$\FI$-object satisfying the hypotheses of \Cref{cogoingdown}, we get analogously :
		\begin{enumerate}
		\item $\delta \bb D \pi_p X \leq \mbf t_p$
		\item $\hmax \bb D \pi_p X \leq \max\{-1, 2\mbf t_p-2, 2\mbf t_{p+1}-2, 2\mbf t_{p+2}-2\}.$
	\end{enumerate}

\end{rmq}
\begin{proof}[Proof of \Cref{cogoingdown}]
Now let us look at what happens for co$\FI$-pointed spaces or spectra. To recover $\FI$-modules in the end we need to dualize and look for the stability of $\bb D \pi_\bullet X = \Hom(\pi_{\bullet} X, A)$ with $A$ an injective abelian group (needed for exactness of the operation of dualizing).  One solution if $X$ is a co$\FI$-spectrum would be to dualize with the Brown-Comenetz spectrum $I_A$ \cite{brown_pontrjagin_1976} that has the property that $\pi_{-k}(\Map(X,I_A)) = \bb D  \pi_k X$.  The dual spectrum $I_A X:= \Map(X, I_A)$ is thus an $\FI$-spectrum and we can apply \Cref{goingdown} and get what we want. This is one solution, but it does not work for pointed spaces as there is no such duality in the unstable homotopy category. Another way is to dualize \emph{after} getting the exact couple inducing the spectral sequence. The hypotheses $\pi_0 = 0$ and $\pi_1$ abelian guarantee that we have a spectral sequence of abelian groups.  For $X$ a co$\FI$-pointed space or spectrum, the cofiltration by the cardinal gives a tower
$$ X \to \cdots \to F_k X \to \cdots \to F_0$$
with fibers the \emph{cofree} objects $M(\bb H X_{\ul k})$. We have a bigraded exact couple $(E_{\bullet,\bullet},A_{\bullet,\bullet} )$ arising from this tower with
$$ E_{p,q} =  \pi_p M(\bb H X)_{\ul q}= M \pi_p (\bb H X)_{\ul q} \text{ and }$$
$$ A_{p,q} = \pi_p F_q X$$
The differential $d_r$ at the $r$-th page is of bidegree $(-1,r)$ with this indexing convention. 
Dualizing every term with respect to the injective abelian group $A$, we get a dual exact couple $\bb D C:= (\bb D E_{\bullet,\bullet} ,\bb D A_{\bullet,\bullet})$ which yields a spectral sequence 
$$E^1_{p,q} = \bb D M(\pi_{-p}(\bb H X_{\ul q})) = M \bb D \pi_{-p}(\bb H X_{ \ul q})$$ because dualizing commutes with finite products. The differential at the $r$-th page is of bidegree $(-1, -r)$. Moreover this spectral sequence converges to $\bb D \pi_{-p}(X)$. To see this, note that by injectivity of $A$, the $r$-th derived exact couple of $\bb D C$ is exactly the dual of the $r$-th derived exact couple of $C$ by induction. The treatment of this spectral sequence is then formally analogous to the one appearing in the $\FI$-version. \\
The last point is checking whether all of this makes sense when $X_S$ and $\bb H X_S$ (once chosen a basepoint of $X_S$) are simply connected for all $S$ yet $X$ has not necessarily a coherent choice of basepoints. 
\begin{rmq}
	This is especially necessary because the co$\FI$ configuration spaces $\conf M$ are not canonically pointed. Actually, they do not even admit a coherent choice of base point at all for $M$ a manifold. This is because there are no homotopy fixed point for the natural action of $S_n$ on $\conf_n(M)$ (see this argument by Connor Malin on mathoverflow: \url{https://mathoverflow.net/a/482985/174025}).
\end{rmq}
To deal with this basepoint problem, we will use some machinery coming from \cite{kupers_representation_2018}.
\begin{defn}\cite[Definition 3.1]{kupers_representation_2018}
	Let $\msf I$ be a small category. 
		\begin{itemize}
			\item[(i)] A \textit{co-$\msf I$-space} \( X \) is a functor \( \msf{I}^{\mathrm{op}} \to \mathbf{Top} \).
			\item[(ii)] A \textit{semistrict co-$\msf I$-space} \( X \) consists of a space \( X_k \) for each object \( k \) of \( \msf{I} \) and for each morphism \( \tau: k \to k' \) a map \( \tau^*: X_{k'} \to X_k \). These must have the property that if \( \sigma: k \to k' \) and \( \tau: k' \to k'' \) are morphisms, then \( \sigma^* \circ \tau^* \) is homotopic to \( (\tau \circ \sigma)^* \) (but these homotopies are not part of the data).
			\item[(iii)] A \textit{homotopy co-$\msf{I}$-space} \( X \) is a functor \( \msf{I}^{\mathrm{op}} \to \mathbf{hTop} \).
		\end{itemize}
	
\end{defn}
These come with their respective notion of map:
\begin{defn}\cite[Definition 3.2]{kupers_representation_2018}
	\begin{itemize}
		\item[(i)] A \textit{map of co-\(\mathsf{I}\)-spaces} \( f: X \to Y \) is a natural transformation of functors \( \mathsf{I}^{\mathrm{op}} \to \mathbf{Top} \).
		
		\item[(ii)] A \textit{semistrict map} \( f: X \to Y \) between semistrict co-\(\mathsf{I}\)-spaces is a sequence of maps \( f_k: X_k \to Y_k \) such that
		\begin{itemize}
			\item[(a)] the equation \( f_k \circ \tau_X^* = \tau_Y^* \circ f_{k'} \) holds, i.e., the following diagram commutes:
			\[
			\begin{tikzcd}
				X_{k'} \arrow[r, "\tau_X^*"] \arrow[d, "f_{k'}"'] & X_k \arrow[d, "f_k"] \\
				Y_{k'} \arrow[r, "\tau_Y^*"'] & Y_k
			\end{tikzcd}
			\]
			
			\item[(b)] for \( \sigma: k \to k' \) and \( \tau: k' \to k'' \), one can choose the homotopies \( \sigma^* \tau^* \sim (\tau \circ \sigma)^* \) for \( X \) and \( Y \) to be compatible with \( f \). More precisely, this condition means that there exist maps 
		\[
		H_{\sigma, \tau}^X: [0,1] \times X_{k''} \to X_k \quad \text{and} \quad H_{\sigma, \tau}^Y: [0,1] \times Y_{k''} \to Y_k
		\]
		\noindent such that the restrictions to \( 0 \) and \( 1 \) give \( \sigma^* \tau^* \) and \( (\tau \circ \sigma)^* \) respectively, and we have
		\[
		f_k \circ H_{\sigma, \tau}^X = H_{\sigma, \tau}^Y \circ (\mathrm{id}_{[0,1]} \times f_{k''});
		\]
		that is, the following diagram commutes:
		\[
		\begin{tikzcd}
			{[0,1] \times X_{k''}} \arrow[r, "H_{\sigma, \tau}^X"] \arrow[d, "\mathrm{id}_{[0,1]} \times f_{k''}"'] & X_k \arrow[d, "f_k"] \\
			{[0,1] \times Y_{k''}} \arrow[r, "H_{\sigma, \tau}^Y"'] & Y_k
		\end{tikzcd}
		\]
		
		\medskip
		
		\noindent (iii) A \textit{map of homotopy co-\(\mathsf{I}\)-spaces} \( f: X \to Y \) is a natural transformation of functors \(\mathsf{I}^{\mathrm{op}} \to \mathbf{hTop} \).

		\end{itemize}
	\end{itemize}
\end{defn}
\begin{prop}\cite[Proposition 3.11]{kupers_representation_2018} \label[prop]{semistrictfiber}
	Let $X \to Y$ a map between semistrict co$\msf I$-spaces. If each \( Y_k \) is \( 1 \)-connected, the collection \( \{ \operatorname{hofib}(f_k, *_{k}) \}_{k} \) can be given the structure of a semistrict co-\(\mathsf{I}\)-space, which we denote by \( \operatorname{hofib}(f) \), such that the map \( p_1: \operatorname{hofib}(f) \to X \) is a semistrict map. This structure is independent of any choice up to objectwise homotopy equivalence.
\end{prop}
With this fact we can deduce the following proposition: 
\begin{lemme}
	Let $X \to Y$ a map between semistrict co-$\msf I$-spaces. Suppose that $X_k$, $Y_k$ and $ D_k = \fiber (X_k \to Y_k)$ are $1$-connected for each $k$. Then the semistrict co-$\msf I$ space structure on the fiber $D$ induces a co-$\msf I$-module structure on $\pi_k D$ and there is a well-defined long exact sequence of co-$\msf I$-modules
	$$ \cdots \to \pi_k D \to \pi_k X \to \pi_k Y \to \pi_{k-1} D \to \cdots$$
\end{lemme}
\begin{proof}
	By \Cref{semistrictfiber} we have a fiber sequence of semistrict $I$-spaces $D \to X \to Y$. As $X$ is simply connected, we can apply a second time \Cref{semistrictfiber} to $D \to X$ and get a fiber sequence $\Omega Y \to D \to X$. Because $D$ is simply connected, applying it a third time yields $\Omega X \to \Omega Y \to D$. Taking homotopy groups of a semistrict $\msf I$-space yields an $\msf I$-abelian group up to conjugacy by $\pi_1$, but the action of $\pi_1$ is trivial on the higher homotopy groups in every space involved (for $X, Y, D$ because $\pi_1$ is trivial and for $\Omega X$ and $\Omega Y$ because an $H$-space has always a trivial action of the $\pi_1$).\\ 
	
	We thus obtain for each $p$ a five-term exact sequence of homotopy groups as $\msf I$-modules 
	$$ \pi_{p+1} X \to \pi_{p+1} Y \to \pi_p D \to \pi_p X \to \pi_p Y$$
	that can assemble and yield the long exact sequence we expect.
\end{proof}
\begin{coroll}\label[coroll]{semistrictSS}
	Let $X \to \cdots X_n \to X_{n-1} \to \cdots X_0$ be a tower of simply connected semistrict $\msf I$-spaces such that each layer $D_n = \fiber (X_n \to X_{n-1})$ is simply connected. Then the spectral sequence of homotopy groups assemble to a spectral sequence of $\msf I$-modules 
	$$ E^1_{p,q} = \pi_p(D_q) \implies \pi_p (X)$$ 
\end{coroll}
\begin{proof}
	The fact that each long exact sequence obtained from the fiber sequences $F_q \to X_q \to X_{q-1}$ lifts to a long exact sequence of $\msf I$-modules implies that we get an bigraded $\msf I$-exact couple $(E_{\bullet, \bullet},A_{\bullet, \bullet})$, with $E_{p,q} = \pi_p D_q$, from which we extract the spectral sequence. 
\end{proof}
Now let us apply this to an unbased, simply connected co$\FI$-space $X$. As $X$ is simply connected, $(F_0 X)_S = X(0)^S$ is simply connected. By assumption, each total fiber $D_q X = F_q X \to F_{q-1} X$ is a simply connected semistrict $\FI$-space which implies by induction that all the $F_q X$ are simply connected. We are thus in the situation where we can apply \Cref{semistrictSS} and get a spectral sequence 
$$ E^1_{p,q} = \pi_p(D_q) \implies \pi_p X.$$
There remains to check that $\pi_p(D_q)$ is still a cofree co$\FI$-module in this unpointed case; we need to do a proof ``by hand'' instead of the categorical argument used in the proof of \Cref{layersfiltration}. Choose a set $S$ and a basepoint $*$ in $F_q X_S$. These propagate to give a basepoint in each $F_q X_T$ for $T \subseteq S$. Write $\bb H Y_U$ for the total fiber of any pointed co-$S$-cube $Y$ at $U \subseteq S$. Let us prove by induction on the cardinal that $$D_q X_S \simeq \prod_{T \subseteq S, |T| = q} \bb H X_T.$$

For every $S$ with $|S| < q$, we have $F_q X_S = *$ and thus $D_q X_S = *$. Then, if $|S| = q$ , we know
    $$ D_q X_S = \tfib_{T \subseteq S} X_T  = X_S$$ by assumption and so $\pi_p F_q X_S = \pi_p \bb H X_{S}$. Finally, if $|S|> k$, both $\bb H F_q X_S$ and $\bb H F_{q-1} X_S$ are points so, as taking total fibers commutes with fibers \cite[Proposition 5.5.4]{munson_cubical_2015}, we get $\bb H D_q X_S =*$ and thus
    $$ D_q X_S \simeq \hocolim_{T \subsetneq S} D_q X_T.$$ But by induction hypothesis $D_qX_T =  \prod_{U \subseteq T, |U| = q} \bb H X_U.$ Therefore we get that $D_q X_S = \prod_{T \subseteq S, |T| = q} \bb H X_T$ and the induction is complete. 
    Applying $\pi_p$ we see that the map $\pi_p D_q X_S\to \pi_p D_qX_T $ is then given by the natural projection $$\prod_{U \subseteq S, |U| = q} \pi_p \bb H X_U \to \prod_{U \subseteq T, |U| = q} \pi_p \bb H X_U $$ and therefore $\pi_p D_q X$ is a cofree co$\FI$-module generated by $\pi_p \bb H X_q$ in cardinal $q$.   

There remains to dualize the spectral sequence with an injective abelian group $A$ to obtain a spectral sequence of $\FI$-modules with
$$ E'^1_{p,q}  = M(\bb D \pi_p{\bb H X_q}) \implies \bb D \pi_p X$$
and run the spectral sequence argument above. This finishes the proof of \Cref{cogoingdown}.
\end{proof}
\begin{rmq}
	We have proved a \gui{going-down} for (co)$\FI$-objects, and a \gui{going-up} that bounds the $\mbf t_k X$s in terms of the bounds on the $\pi_k X$, the latter working only for globally bounded-below $\FI$-spectra. The author has not been able to generalize this to general $\FI$-spectra. For example, we don't know if having $t_i \pi_{-p} X \leq ap+b$ for $i=0,1$ implies that the $\mbf t_{-p} X$ are finite, and if there is an explicit bound.  
\end{rmq}
\begin{rmq}
	Note that the going-down theorem implies that if $X$ is a stable co$\FI$-spectrum and $E$ a coconnective spectrum, then the $E$-cohomology groups of $X$ are representation stable, because they are the homotopy groups of the $\FI$-mapping spectrum $\mathrm{Map}(X,E)$. As $\bb H \mathrm{Map}(X,E) = \mathrm{Map}(\bb H X,E)$ the mapping spectrum is stable and the going down theorem gives explicit bounds on the stability of its homotopy groups. If $E$ is a connective spectrum, this argument does not work but in this case the Atiyah-Hirzebruch spectral sequence is a first-quadrant spectral sequence and we can control the generation degrees with a type A spectral sequence argument from \cite{church_linear_2018}.
\end{rmq}

\section{Representation stability and Blakers-Massey theorem}\label{Applications}

In this section we apply the previous general results to the co$\FI$-space of configurations on a manifold and its \gui{$\FI$-cohomology}, in other words the data of the homotopy fibers of all the cubes we can define by forgetting points. We will see that we can deduce quantitative bounds on the $\pi_k \conf M$ from the study of these cubes and the generalized Blakers-Massey theorem.
\subsection{Cartesianity of the cube of configuration spaces}
This approach begins with the following proposition coming from \cite[Example 6.2.9.]{munson_cubical_2015}, 
\begin{prop}\label[prop]{blakersmasseyconfigurations}
    Let $M$ be an $u$-connected $d$-dimensional topological manifold. Consider the $n$-cube $C = (S \mapsto \conf_{\ul n - S}(M))$. Then the cube $C$ is $((n-1)(d-2) + 1)$-cartesian.
\end{prop}
\begin{proof}
This fact is stated for smooth manifolds in \cite{munson_cubical_2015} but the assumption is actually not needed. The only point to check in the proof is the fact that for a topological $d$-manifold $M$ and a point $x\in M$, the inclusion $M-\{x\} \subseteq M$ is still $(d-1)$-connected. This is true by \cite[Lemma 3.3]{GGLG_2017}.
\end{proof}
We now recall below the most general statement of the generalized Blakers-Massey theorem for $n$-cubes. Following \cite{munson_cubical_2015}, for an $n$-cube $X_\bullet$ and a subset $U \subset \ul n$, we write $\partial^U X$ for the subcube $S \mapsto X_S$ for $S \subseteq U$ and $\partial_U X$ for the subcube $S \mapsto X_S$ for $S \supseteq \ul n \backslash U$.  
\begin{thm}\cite[Theorem 2.5]{goodwillie_calculus_2003}\label{BMcart}
	Let $\mathcal{X} = (T \mapsto X_T)$ be an $S$-cube with $|S| = n \geq 1$. Suppose that
	\begin{enumerate}
		\item for each $\vide \neq U \subset S$ the $U$-cube $\partial^U \mathcal{X}$ is $k_U$-cocartesian, and
		\item for $U \subset V$, $k_U \leq k_V$.
	\end{enumerate}
	Then $\mathcal{X}$ is 
	\[
	(1 - n + \min\{\Sigma_\alpha k_{T_\alpha}: \{T_\alpha\} \text{ is a partition of } S\})\text{-cartesian}.
	\]
	In particular,  if $X$ is a strongly cocartesian cube, then $X$ is $(1-n + \sum_{i \in \ul n} k_i)$-cartesian, writing $k_i$ for the connectivity of $X_\vide \to X_{\{i\}}$.
\end{thm}
Dually, 
\begin{thm}\cite[Theorem 2.6]{goodwillie_calculus_2003}\label{BMcocart}
	Let $\mathcal{X} = (T \mapsto X_T)$ be an $S$-cube with $|S| = n \geq 1$. Suppose that
	\begin{enumerate}
		\item for each $\vide \neq U \subset S$, the $U$-cube $\partial_{S-U} \mathcal{X}$ is $k_U$-cartesian, and
		\item for $U \subset V$, $k_U \leq k_V$.
	\end{enumerate}
	Then $\mathcal{X}$ is 
	\[
	(-1 + n + \min\{\Sigma_\alpha k_{T_\alpha}: \{T_\alpha\} \text{ is a partition of } S\})\text{-cocartesian}.
	\]
	In particular,  if $X$ is a strongly cartesian cube, then $X$ is $(-1+n + \sum_{i \in \ul n} k_i)$-cocartesian, writing $k_i$ for the connectivity of $X_{\ul n} \to X_{\ul n \setminus \{i\}}$.
\end{thm}
We give now a "cartesian ranges to $\FI$-ranges" result that is an intermediate step to prove \Cref{thmA} and \Cref{thmB}. 
\begin{defn}
	For $f : \bb N \to \bb N$, we define $\mc P f : \bb N \to \bb N$ by $$\mc Pf (n) =  \min \left\{ \sum_{T \in \alpha} f(|T|) : \alpha \in \msf{Par}(\ul n)\right\} $$ where $\msf{Par}(\ul n)$ is the set of partitions of the finite set $\ul n$. 
\end{defn}
\begin{lemme}\label[lemma]{cartesiantoFI}
	Let $X$ be a connected co$\FI$-space. If $f: \bb N\backslash \{0\} \to \bb N$ is a weakly increasing function such that the $n$-cube 
	$$ \ul n  \supset T \mapsto X_{\ul n \backslash T}$$ is $f(n)$-cartesian, then the following hold : 
	\begin{enumerate}
		\item Suppose $\pi_1 X_S$ is abelian for all $S$. Then for every $p$,
		$$ \deg \pi_p \bb HX \leq \sup (\{n \text{ such that } f(n) \leq p\}\cup \{-1\})$$
  
		\item For every $p$, 
		$$ \mbf{ t}_{-p} (C^{-*} X) \leq \sup (\{n : \mc Pf(n) \leq p\}\cup \{-1\})$$ 
	\end{enumerate}
\end{lemme}
\begin{proof}
	Part $1$ is elementary : by assumption, $\bb H X_n$ is $f(n)-1$-connected as it is the total fiber of the $n$-cube $T \mapsto X_{\ul n \backslash T }$. This means $\pi_p \bb H X_n = 0$ for $p \leq f(n)-1$ i.e the biggest $n$ for which $\pi_p \bb H X_n \neq 0$ is $\max \{ n : f(n) \leq p\}$.  For part $2$, we know by \Cref{BMcocart} that the $n$-cube $T \mapsto X_{\ul n \bs T}$ is $(1-n+\mc Pf (n))$-cocartesian. Applying the singular cochain functor, which sends homotopy colimits to homotopy limits, and preserves connectivity, we get a $(1-n+\mc Pf(n))$-cartesian cube of cochain complexes, and hence a $(1+\mc Pf(n))$-cocartesian cube, because the total fiber is the $n$-th desuspension of the total cofiber in a stable category. We then apply the same reasoning as part $1)$.
\end{proof}
Let us now prove \Cref{thmA} when the manifold $M$ is simply connected. Using the language of $\FI$-homology, this can be seen as a result about the slope of the co$\FI$-homology $\bb H \conf M$. Applying \Cref{cartesiantoFI}, part $1$ we get : 
$$ \mbf t_p \leq \max \{ n \text{ such that } (n-1)(d-2) < p\} = \ceil* {\frac{p}{d-2}} $$

We can thus apply the going-down \Cref{cogoingdown}, spatial, unbased version to the co$\FI$-space $X:= \conf M$  where we take $M$ to be a simply-connected, $d$-dimensional manifold with $d \geq 3$. The connectivity hypotheses for the \Cref{semistrictSS} are satisfied because $X_S = \conf_S M$ is simply connected for each $S$ and the total fiber of configuration spaces is $(d-2)(|S|-1)$-connected by \Cref{blakersmasseyconfigurations}. So for $S$ of cardinal $\geq 2$, and $d \geq 3$, the fiber is simply-connected. For $|S| = 1$, note that the fiber is $\conf_1 M = M$ which is simply connected.  Applying \Cref{rmq2.8}, we obtain for any injective abelian group $A$, (writing $\bb D V$ for $\Hom(V,A)$) the bounds
$$ \delta \bb D \pi_p X \leq \floor*{\frac{p-1}{d-2}}+1$$
$$ \hmax \bb D \pi_p X \leq 2 \floor*{\frac{p-1}{d-2}}$$

But we would like to get bounds for $\Hom(\pi, \bb Z)$ and $\Ext(\pi, \bbZ)$ as in \cite{kupers_representation_2018}. This is possible exploiting the exact sequence 
$$ 0 \to \Hom(\pi, \bb Z) \to \Hom(\pi, \bb Q) \to \Hom(\pi, \bb Q/\bb Z) \to \Ext(\pi, \bb Z) \to 0 $$
which is functorial in the abelian group $\pi$. 
We now use \Cref{localmaxexact} to obtain:
\begin{enumerate}
    \item{$\delta \Hom(\pi_p X, \bb Z) \leq \delta \Hom(\pi_p X, \bb Q) \leq  \floor* {\frac{p-1}{d-2}}+1 $}
    \item{$\hmax \Hom(\pi_p X, \bb Z) \leq  2\floor*{ \frac{p-1}{d-2}} $}
    \item{$\delta \Ext(\pi_p X, \bb Z) \leq \delta \Hom(\pi_p X, \bb Q/ \bb Z) \leq  \floor*{\frac{p-1}{d-2}}+1 $}
    \item{$\hmax \Ext(\pi_p X, \bb Z) \leq  2 \floor*{\frac{p-1}{d-2}}. $}
\end{enumerate}
Combining this with \Cref{bahranlemma}, we complete the proof of \Cref{thmA} in the simply connected case. Note that we prefer to write the main result in terms of floor function instead of ceil functions. The conversion is easy because of the following elementary lemma. 
\begin{lemme}\label[lemme]{floorceil}
	For $a,b\in \bb Z$, we have
	$$ \ceil* {\frac {a+1} b} = \floor*{\frac a b}+1.$$
\end{lemme}
\subsection{Recovering stability of the cohomology of configuration spaces}
Another interesting consequence of \Cref{blakersmasseyconfigurations} is that we can recover a bound for the generation and presentation degrees of $H^*(\conf M, \bb Z)$ applying the Blakers-Massey theorem a second time. More concretely we obtain the following theorem.
\begin{thm*}[\Cref{thmB}]
    Let $M$ be an $u$-connected ($u\geq 0$) manifold of dimension $d \geq 3$. Then its integral cohomology groups are stable with the following bounds:
\begin{enumerate}
	\item{$t_0 H^p \conf M \leq \begin{cases} 
			2\floor*{\frac{2p}{d-1}} , & \text{if } u+1 \geq \frac{d-1}{2}, \\[6pt]
			2\floor*{\frac{p}{u+1}} & \text{otherwise.} \\[6pt]
		\end{cases}$}
	\item{$t_1 H^p \conf M \leq \begin{cases} 
			2\floor*{\frac{2p}{d-1}}+1, & \text{if } u+1 \geq \frac{d-1}{2}, \\[6pt]
			2\floor*{\frac{p}{u+1}}+1 & \text{otherwise.} \\
		\end{cases}$}
\end{enumerate}
\end{thm*}
\begin{proof}
	By \Cref{blakersmasseyconfigurations} we know that for any finite set $S$ the $S$-cube 
	$$ S \supseteq T \mapsto \conf_{S \setminus U} M $$ is $(|T|-1)(d-2)+1$-cartesian. Moreover as $M$ is $u$-connected, the map $M = \conf_{1} M \to \star$ is $(u+1)$-connected. This means that for every $T \subset S$, the sub-cube $\partial_T \conf_\bullet M$ is $k_T$ cartesian where:
	\[
	k_T =
	\begin{cases} 
		(|T|-1)(d-2) + 1, & \text{if } |T| \geq 2, \\
		u+1 & \text{if } |T| = 1. \\
	\end{cases}
	\]
Now we apply \Cref{cartesiantoFI}, part $2$, with the function $f$ given by 
\[
f(n)=
\begin{cases}
	(n-1)(d-2)+1, & n\ge 2,\\[4pt]
	u+1, & n=1.
\end{cases}
\]

Fix $T_\alpha$ a partition of $\ul n$. Let $s$ be the number of singletons in this partition and $t$ the total number of elements in the partition. We have by basic combinatorics $s + 2(t-s) \leq n$ so $t-s \leq \frac{n-s}{2}$. Moreover we get 
\begin{align*}
	 \sum_{\alpha} k_{T_\alpha} &= s (u+1) + \sum_{\alpha \text{ s.t. } |T_\alpha| \geq 2} ((|T|-1)(d-2) + 1) \\ &= s(u+1)  + (n-s)(d-2) - (t-s)(d-3) \\
	&= n(d-2) +  s((u+1) - d+ 2) - (t-s) (d-3) \\
	&\geq n(d-2) + s(u+1-d+2) - \frac{n-s}{2} (d-3) \\
	&= n \frac{d-1}{2} + s \left(u+1 - \frac{d-1}{2}\right)
\end{align*}
Taking the minimum of this value for $s$ ranging from $0$ to $n$ gives then 
	\[ \mc P f(n) = 
\min\{\Sigma_\alpha k_{T_\alpha}: \{T_\alpha\} \text{ is a partition of } \ul n \} \geq 	
\begin{cases} 
	n \frac{d-1}{2}, & \text{if } u+1 \geq \frac{d-1}{2}, \\
	n(u+1) & \text{otherwise.} \\
\end{cases}
\]
Applying \Cref{cartesiantoFI}, part $2$, we get that the singular cochain complex (seen as a negatively-graded chain complex) is an $\FI$-complex stable with $$\mbf t_{-p} C^{-*}(\conf M, \bb Z) \leq  \begin{cases} 
	\floor*{\frac{2p}{d-1}}, & \text{if } u+1 \geq \frac{d-1}{2}, \\[6pt]
	\floor*{\frac{p}{u+1}} & \text{otherwise.} \\
\end{cases}$$
Running the going-down theorem, stable version, we get the expected bounds of
\begin{enumerate}
	\item{$t_0 H^p \conf M \leq \begin{cases} 
			2\floor*{\frac{2p}{d-1}}  , & \text{if } u+1 \geq \frac{d-1}{2}, \\[6pt]
			2\floor*{\frac{p}{u+1}}  & \text{otherwise.} \\[6pt]
		\end{cases}$}
	\item{$t_1 H^p \conf M \leq \begin{cases} 
			2\floor*{\frac{2p}{d-1}}+1, & \text{if } u+1 \geq \frac{d-1}{2}, \\[6pt]
		2\floor*{\frac{p}{u+1}}+1 & \text{otherwise.} \\
		\end{cases}$}

\end{enumerate}

\end{proof}
Note that the slope is inversely proportional to the connectivity (or the dimension depending on the cases) of the manifold we consider. This phenomenon, and the disjunction at half the dimension has already been observed in Bahran's recent paper \cite{bahran_regularity_2024}.  Note that for oriented manifolds, by Poincaré duality the first case $u \geq \frac{d-1}{2}$ happens only if $M$ is contractible or a homotopy $d$-sphere (and hence a sphere by Poincaré conjecture). 
\begin{rmq}
	Note that in the case $d=2$, we cannot run the Blakers-Massey argument anymore because in this case the minimum value of $\sum_\alpha k_{T_\alpha} $ is $1$. Given this we can only deduce that the $n$-cube $S \mapsto C_*(\conf_S M, \bb Z)$ is $1$-cartesian which is not sufficient to run the going-down theorem and get ranges on the cohomology. However, Miller and Wilson dealt with the $2$-dimensional case and provided ranges for the $\FI$-hyperhomology of $C^*(\conf M, \bbZ)$ for $M$ a $2$-dimensional manifold, see the proof of \cite[Theorem 1.1]{miller_fihyperhomology_2019}.
\end{rmq}
\subsection{$\FI_G$-modules and non-trivial fundamental groups}
In this subsection we show how the previous constructions and arguments generalize when replacing $\FI$ with the category $\FI_G$ of $G$-sets with finitely many orbits and injective $G$-equivariant maps, for $G$ a discrete group. As for $\FI$-modules, there is a well-defined ``$\FI_G$-homology`` 
$$ \bb H^{\FI_G}: \mbf{Fun}(\FI_G, \mc C) \to \mbf{Fun}(\FB_G, \mc C)$$ defined as the derived left adjoint of the extension-by-zero functor. Here $\mc C$ denotes a pointed $\infty$-category. This $\FI_G$-homology was studied extensively in \cite{sam_representations_2019}.  The main results of this section are the following. 
\begin{prop}\label[prop]{prop:FIGeqFI}
	Let $X$ be a $\FI_G$-object in a pointed $\infty$-category $\mc C$. Then the $\FI_G$-homology $\bb H^{\FI_G}$ of $X$ is equivalent to the $\FI$-homology $\bb H^{\FI}$ of $X$, seen as an $\FI$-object through the inclusion $\FI \inj \FI_G$ that sends a finite set $S$ to the free $G$-set $S \times G$. 
\end{prop}
Consider $\FI \xinj{i} \FI_G$ and $\FB \xinj{j} \FB_G$ the inclusions. We want to show that the following square commutes. 
\[\begin{tikzcd}
	{\mbf{Fun}(\FI_G, \mc C)} & {\mbf{Fun}(\FB_G, \mc C)} \\
	{\mbf{Fun}(\FI, \mc C)} & {\mbf{Fun}(\FB, \mc C)}
	\arrow["{\bb H^{\FI_G}}", from=1-1, to=1-2]
	\arrow["{i^*}"', from=1-1, to=2-1]
	\arrow["{j^*}", from=1-2, to=2-2]
	\arrow["{\bb H^{\FI}}"', from=2-1, to=2-2]
\end{tikzcd}\]
It suffices to show that the square of right adjoints commutes: 
\[\begin{tikzcd}
	{\mbf{Fun}(\FI_G, \mc C)} & {\mbf{Fun}(\FB_G, \mc C)} \\
	{\mbf{Fun}(\FI, \mc C)} & {\mbf{Fun}(\FB, \mc C)}
	\arrow["\beta^*"', from=1-2, to=1-1]
	\arrow["{\mathrm{Ran}_i}", from=2-1, to=1-1]
	\arrow["{\mathrm{Ran}_j}", from=2-2, to=1-2]
	\arrow["\alpha^*"', from=2-2, to=2-1]
\end{tikzcd}\]
Where $\alpha^*$ and $\beta^*$ are the corresponding extension-by-zero functors. 
Let us compute the right Kan extension of $i$ and $j$. 
Take $F \in \mbf{Fun}(\FI, \mc C)$. Then 
$$ \mathrm{Ran}_i F(S \times G) = \holim_{S \times G \downarrow \FI} F$$
Let us study the comma category $\mc D = (S \times G) \downarrow \FI$. An object is given by a pair $(T,u)$ where $T$ is a finite set and $u$ is a $G$-injection $S \times G \to T \times G$. A morphism $(T, u) \to (T',u')$ is an injection of finite sets $T \to T'$ such that the obvious diagram commutes. Let us show that $\mc D$ is a disjoint union of categories admitting an initial object. First of all notice that a $G$-injection  $u: S \times G \to T \times G$ induces a map $\hat u: S \to G$ by $\hat u(s) =p_G(u(s,1_G))$ with $p_G: T \times G \to G$ the projection. 
Define the full subcategory $\mc D^f$ of $\mc D$ consisting of objects $(T,u)$ such that $\hat u = f$. The collection $(\mc D^f)_{f \in G^S}$ partitions the objects of $\mc D$. Now notice that if there is a map $(T,u) \to (T',u')$ in $\mc D$ their functions $\hat u$ and $\hat u'$ are equal. Therefore there are no maps between objects of $\mc D^f$ and of $\mc D^g$ if $f \neq g$. Finally, notice that each of the $\mc D^f$ has an initial object given by $$S \times G \to S\times G$$  
$$ (s,g) \mapsto (s, g\cdot f(s))$$

From this argument we deduce that $\mathrm{Ran}_i F(S \times G) = \holim_{S \times G \downarrow \FI} F =  F(S\times G)^{G^S}$. The same argument works for $\FB \inj \FB_G$ and we get $\mathrm{Ran}_j F( S \times G) = F(S \times G)^{G^S}$ where $F$ is an $\FB$-object and $S$ an arbitrary finite set. From this description we see easily that the above square of right adjoints commutes. 
\begin{defn}
	Let $X$ be a free $G$-space with $G$ a discrete group. The \emph{orbit configuration space} $\conf^G_n(X)$ of $n$ points in $X$ is the space of configurations of $n$ points with disjoint $G$-orbits. In other words, it is the space 
	$$\mathrm{Emb}^G(\ul n \times G, X)$$ of $G$-equivariant embeddings from the discrete $G$-space $\ul n \times G$ into $X$. 
\end{defn}
\begin{lemme}\label[lemme]{lem:bmG}
	Let $G$ be a discrete group and $M$ be a $d$-dimensional manifold with a free and properly discontinuous $G$-action. Then the $n$-cube of orbit configurations
	$$\ul n \supseteq S \mapsto \conf^G_{\ul n- S}(M)$$
	is $((n-1)(d-2)+1)$-cartesian. 
\end{lemme}
\begin{proof}
	We mimic the proof of \cite[Example 6.2.9]{munson_cubical_2015} for the ordinary configuration spaces. 	See the $n$-cube $S \mapsto \conf^G_{\ul n - S}(M)$ as a map of $(n-1)$-cubes: 
		$$ (R \mapsto \conf^G_{\{n\} \cup \ul{n-1} \backslash R}(M)) \to (R \mapsto \conf^G_{ \ul{n-1} \backslash R}(M))$$ 
		where $R \subseteq \ul{n-1}$. Choose a basepoint $(x_1, \ldots, x_{n-1}) \in \mathrm{Conf}^G_{n-1}(M)$. This gives rise to a basepoint in $\mathrm{Conf}^G_{\ul {n-1}-R}(M)$ for all $R \subseteq \ul {n-1}$.  The map of $(n-1)$-cubes, as an $n$-cube, is as cartesian as the $(n-1)$-cube of homotopy fibers.  By work of Xicot{\'e}ncatl \cite[Theorem 2.2]{xico}, for each $R$, the restriction map $\mathrm{Conf}^G_{\{n\}\cup \ul {n-1} \backslash R}(M) \to \mathrm{Conf}^G_{\ul {n-1} \backslash R}(M) $ above is a fibration, whose fiber is $M \bs \bigsqcup_{i \in \ul {n-1} \backslash R} Gx_i $, which we write as $M \bs (G\times(n-1-R))$ for short. The $(n-1)$-cube of fibers is
		\[
		R \mapsto M\backslash (G \times (n-1-R)),
		\]
		with $M \bs (G \times (n-1-R)) \to M \bs (G \times (n-1-S))$ the evident inclusion for $R \subseteq S$. Further, $M\bs (G \times (n-1-R))$ is an open subset of $M \bs (G \times (n-1-S))$ since it is the complement of the closed subset $S-R$. Each square face of this $(n-1)$-cube is of the form
		\[
		\begin{tikzcd}
				M \bs (G \times (n-1-R)) \arrow[r] \arrow[d] & M \bs (G\times (n-1-(R \cup \{i\}))) \arrow[d] \\
				M\bs (G \times (n-1-(R \cup \{j\}))) \arrow[r] & M \bs (G \times (n-1-(R \cup \{i, j\})))
			\end{tikzcd}
		\]
		for $R \subset n-1$, $i, j \notin R$, and $i \neq j$. It is clear by inspection that
		\[
		M\bs (G \times (n-1-R)) = M \bs (G\times(n-1-(R \cup \{i\}))) \cap M \bs (G \times (n-1-(R \cup \{j\})))
		\]
		and
		\[
		M \bs ((n-1-G \times (R \cup \{i\}))) \cup M\bs (G \times (n-1-(R \cup \{j\}))) = M \bs (G \times (n-1-(R \cup \{i, j\}))).
		\]
		Hence this square is homotopy cocartesian, and so the $(n-1)$-cube in question is strongly homotopy cocartesian. The maps
		\[
		M \bs (G \times (n-1)) \to M \bs( G \times (n-1-\{i\}))
		\]
		are $(d-1)$-connected for each $i \in \ul n-1$, because the cofiber is a wedge of $d$-spheres. By the Blakers-Massey theorem, this $(n-1)$-cube, and hence the original $n$-cube, is $((n-1)(d-2)+1)$-cartesian.
\end{proof}
\begin{thm}\label{thm:Gequivariantconfigurations}
	Let $G$ be a discrete group and $M$ be a simply connected manifold of dimension $d \geq 3$ with a properly discontinuous free $G$-action. Then the collection of orbit configurations $\conf^G_\bullet M$ has the structure of a co$\FI_G$-space and its homotopy groups are stable with the same bounds as in \Cref{thmA}. 
\end{thm}
\begin{proof}
	This results from the combination of \cref{lem:bmG} and \cref{prop:FIGeqFI} together with the going-down theorem for co$\FI$-objects (\cref{cogoingdown}). What makes the argument work is that $\conf^G_\bullet M$ is simply connected. The computations are the same as for the $\FI$-case. 
\end{proof}
\begin{coroll}
	Let $M$ be a $d$-dimensional connected manifold with fundamental group $G$ and $d \geq 3$. Then the homotopy groups $\pi_p \conf M$ for $p \geq 2$ can be given the structure of a co$\FI_G$-module, and their duals have finite generation degree with the same bounds as in $\Cref{thmA}$. 
\end{coroll}
\begin{proof}
 The universal covering of $\conf_n M$ is $\conf^G_n(\widehat M)$ where $\widehat M$ is the universal covering of $M$ endowed with its $G$-action. See for example \cite[Proposition 5.6]{kupers_representation_2018}. The action of $G^n$ on it is properly discontinuous and free, for example see \cite[Proposition 1.40]{hatcher} for a proof. We can then apply \Cref{thm:Gequivariantconfigurations} and use the fact that the higher homotopy groups of the universal covering are the same as the ones of the base. 
\end{proof}
\begin{rmq}
	We can also recover finite generation degree for the integral cohomology of orbit configuration spaces with explicit ranges, combining \cref{lem:bmG} with a second application of the Blakers-Massey theorem as in the proof of \Cref{thmB}. The ranges we get are exactly the same as in \Cref{thmB}.
\end{rmq}
\bibliographystyle{alpha}
\bibliography{bib_fispectra.bib}
\end{document}